\documentclass{llncs}
\usepackage{makeidx}  % allows for indexgeneration
\usepackage{amsmath}
\usepackage{amsfonts,amssymb}
\usepackage[thinlines,thiklines]{easybmat}
\usepackage{multirow}
\usepackage{algorithmic}
\usepackage{algorithm}

\newcommand{\Z}{\mathbb{Z}} % les entiers relatifs
\newcommand{\p}{\mathfrak{p}} % un p gothik
\renewcommand{\algorithmicrequire}{\textbf{Input:}}
\renewcommand{\algorithmicensure}{\textbf{Output:}}
\begin{document}

\mainmatter              % start of the contributions
\title{Practical improvements to class group and regulator
computation of real quadratic fields}
\titlerunning{Improvements in real quadratic number fields}  % abbreviated title (for running head)
\author{Jean-Fran\c{c}ois Biasse\inst{1} \and Michael J. Jacobson,
Jr.\inst{2}\thanks{The second author is supported in part by NSERC of
Canada.}}
\authorrunning{J-F Biasse and M. J. Jacobson, Jr.}   % abbreviated author list (for running head)
\tocauthor{Jean-Fran\c{c}ois Biasse  Michael J. Jacobson, Jr.}
\institute{ \'{E}cole Polytechnique, 91128 Palaiseau, France\\
\email{biasse@lix.polytechnique.fr}
\and
Department of Computer Science, University of Calgary \\
2500 University Drive NW, Calgary, Alberta, Canada T2N 1N4\\
\email{jacobs@cpsc.ucalgary.ca}}

\maketitle              % typeset the title of the contribution

\begin{abstract}
We present improvements to the index-calculus algorithm for the
computation of the ideal class group and regulator of a real quadratic
field. Our improvements consist of applying the double large prime
strategy, an improved structured Gaussian elimination strategy, and
the use of Bernstein's batch smoothness algorithm. We achieve a
significant speed-up and are able to compute the ideal class group
structure and the regulator corresponding to a number field with a
110-decimal digit discriminant.
\end{abstract}

\section{Introduction}
Computing invariants of real quadratic fields, in particular the ideal
class group and the regulator, has been of interest since the time of
Gauss, and today has a variety of applications.  For example, solving
the well-known Pell equation is intimately linked to computing the
regulator, and integer factorization algorithms have been developed
that make use of this invariant.  Public-key cryptosystems have also
been developed whose security is related to the presumed difficulty of
these computational tasks.  See \cite{JWPellBook} for details.

The fastest algorithm for computing the ideal class group and
regulator in practice is a variation of Buchmann's index-calculus
algorithm \cite{BSub} due to Jacobson \cite{JacobsonPhd}.  The
algorithm on which it is based has subexponential complexity in the
size of the discriminant of the field.  The version in
\cite{JacobsonPhd} includes several practical enhancements, including
the use of self-initialized sieving to generate relations, a single
large-prime variant (based on that of Buchmann and D\"ullman
\cite{BDDist} in the case of imaginary quadratic fields), and a
practical version of the required linear algebra.  This approach
proved to work well, enabling the computation of the ideal class group
and regulator of a real quadratic field with a 101-decimal digit
discriminant \cite{JSWkeyex}.  Unfortunately, both the complexity
results of Buchmann's algorithm and the correctness of the output are
dependent on the Generalized Riemann Hypothesis (GRH).  Nevertheless,
for fields with large discriminants, this approach is the only one
that works.

Recently, Biasse \cite{biasse} presented practical improvements to the
corresponding algorithm for imaginary quadratic fields.  These
included a double large prime variant and improved algorithms for the
required linear algebra.  The resulting algorithm was indeed faster
then the previous state-of-the-art \cite{JacobsonPhd}, and enabled the
computation of the ideal class group of an imaginary quadratic field
with $110$ decimal digit discriminant.

In this paper, we describe a number of practical improvements to the
index-calculus algorithm for computing the class group and regulator
of a real quadratic field.  In addition to adaptations of Biasse's
improvements in the imaginary case, we have found some modifications
designed to improve the regulator computation part of the algorithm.
We also investigate applying an idea of Bernstein \cite{bernstein} to
factor residues produced by the sieve using a batch smoothness test.
Extensive computations demonstrating the effectiveness of our
improvements are presented, including the computation of class group
and regulator of a real quadratic field with $110$ decimal digit
discriminant.

This paper is organized as follows.  In the next section, we briefly
recall the required background of real quadratic fields, and give an
overview of the index-calculus algorithm using self-initialized
sieving.  Our improvements to the algorithm are described in
Section~\ref{sec:improvements}, followed by numerical results in
Section~\ref{numerical}.

\section{Real Quadratic Fields}

We present an overview of required concepts related to real quadratic
fields and the index-calculus algorithm for computing invariants.  For
more details, see \cite{JWPellBook}.

Let $K = \mathbb{Q}(\sqrt{\Delta})$ be the real quadratic field of
discriminant $\Delta,$ where $\Delta$ is a positive integer congruent
to $0$ or $1$ modulo $4$ with $\Delta$ or $\Delta/4$ square-free.  The
integral closure of $\mathbb{Z}$ in $K$, called the maximal order, is
denoted by $\mathcal{O}_{\Delta}.$ An interesting aspect of real quadratic
fields is that their maximal orders contain infinitely many
non-trivial units, i.e., units that are not roots of unity.  More
precisely, the unit group of $\mathcal{O}_\Delta$ consists of an order $2$
torsion subgroup and an infinite cyclic group.  The smallest unit
greater than $1,$ denoted by $\varepsilon_\Delta,$ is called the
fundamental unit.  The regulator of $\mathcal{O}_\Delta$ is defined as
$R_\Delta = \log \varepsilon_\Delta.$

The fractional ideals of $K$ play an important role in the
index-calculus algorithm described in this paper.  In our setting, a
fractional ideal is a rank $2$ $\mathbb{Z}$-submodule of $K.$  Any
fractional ideal can be represented as
\begin{equation*}
\mathfrak{a} = \frac{s}{d} \left[a \mathbb{Z} + \frac{b +
    \sqrt{\Delta}}{2} \mathbb{Z} \right] \enspace,
\end{equation*}
where $a,b,s,d \in \Z$ and $4a ~|~ b^2 - \Delta.$ The integers
$a,$ $s,$ and $d$ are unique, and $b$ is defined modulo $2a.$ The
ideal $\mathfrak{a}$ is said to be primitive if $s = 1,$ and $d
\mathfrak{a} \subseteq \mathcal{O}_\Delta$ is integral.  The norm of
$\mathfrak{a}$ is given by $\mathcal{N}(\mathfrak{a}) = a s^2 / d^2.$

Ideals can be multiplied using Gauss's composition formulas for
indefinite binary quadratic forms.  Ideal norm respects ideal
multiplication, and the set $\mathcal{I}_\Delta$ forms an infinite abelian
group with identity $\mathcal{O}_\Delta$ under this operation.  The
inverse of $\mathfrak{a}$ is
\begin{equation*}
\mathfrak{a}^{-1} = \frac{d}{s a} \left[a \mathbb{Z} + \frac{-b +
    \sqrt{\Delta}}{2} \mathbb{Z} \right] \enspace.
\end{equation*}
The group $\mathcal{I}_\Delta$ is generated by the prime ideals of
$\mathcal{O}_\Delta,$ namely those integral ideals of the form $p
\mathbb{Z} + (b_p + \sqrt{\Delta})/2 \mathbb{Z}$ where $p$ is a prime
that is split or ramified in $K.$ As $\mathcal{O}_\Delta$ is a
Dedekind domain, the integral part of any fractional ideal can be
factored uniquely as a product of prime ideals.  To factor
$\mathfrak{a},$ it suffices to factor $\mathcal{N}(\mathfrak{a})$ and,
for each prime $p$ dividing the norm, determine whether the prime
ideal $\mathfrak{p}$ or $\mathfrak{p}^{-1}$ divides $\mathfrak{a}$
according to whether $b \equiv b_p$ or $-b_p$ modulo $2p.$

The ideal class group, denoted by $Cl_\Delta,$ is the factor group
$\mathcal{I}_\Delta / \mathcal{P}_\Delta,$ where $\mathcal{P}_\Delta
\subseteq \mathcal{I}_\Delta$ is the subgroup of principal ideals.
The class group is finite abelian, and its order is called the class
number, denoted by $h_\Delta.$ By computing the class group we mean
computing the elementary divisors $m_1, \dots, m_l$ with $m_{i+1} ~|~
m_i$ for $1 \leq i < l$ such that $Cl_\Delta \cong \Z/m_1 \Z \times
\dots \times \Z / m_l \Z.$

\subsection{The Index-Calculus Algorithm}\label{sieving}

Like other index-calculus algorithms, the algorithm for computing the
class group and regulator relies on finding certain smooth quantities,
those whose prime divisors are all small in some sense.  In the case
of quadratic fields, one searches for smooth principal ideals for
which all prime ideal divisors have norm less than a given bound
$B_1.$ The set of prime ideals $\mathcal{B} = \{ \mathfrak{p}_1,
\dots, \mathfrak{p}_n \}$ with $\mathfrak{p}_i \leq B_1$ is called the
factor base.

A principal ideal $(\alpha) = \mathfrak{p}_1^{e_1} \dots
\mathfrak{p}_n^{e_n}$ with $\alpha \in K$ that factors completely over
the factor base yields the relation $(e_1,\dots,e_n, \log |\alpha|).$
The key to the index-calculus algorithm is the fact, proved by
Buchmann \cite{BSub}, that the set of all relations forms a sublattice
$\Lambda \subset \mathbb{Z}^n \times \mathbb{R}$ of determinant
$h_\Delta R_\Delta$ provided that the prime ideals in the factor base
generate $Cl_\Delta.$ This follows, in part, due to the fact that $L,$
the integer component of $\Lambda,$ is the kernel of the homomorphism
from $\mathbb{Z}^n$ to $Cl_\Delta$ given by $\mathfrak{p}_1^{e_1}
\dots \mathfrak{p}_n^{e_n}$ for $(e_1,\dots,e_n) \in \mathbb{Z}^n.$ If
$\mathfrak{p}_1, \dots, \mathfrak{p}_n$ generate $Cl_\Delta,$ then
this homomorphism is surjective, and the homomorphism theorem then
implies that $\mathbb{Z}^n / L \cong Cl_\Delta.$

The main idea behind the index-calculus algorithm is to find random
relations until they generate the entire relation lattice $\Lambda.$
Let $\Lambda'$ denote the sublattice of $\Lambda$ generated by the
relations that have been computed.  To determine whether $\Lambda' =
\Lambda,$ one computes an approximation $h^*$ of $h_\Delta R_\Delta$
such that $h^* < h_\Delta R_\Delta < 2 h^*.$ The value $h^*$ is
obtained by approximating the $L$-function $L(1,\chi_\Delta),$ where
$\chi_\Delta$ denotes the Kronecker symbol $(\Delta/p),$ and applying
the analytic class number formula.  If $\Lambda' \subset \Lambda,$
then $\det (\Lambda')$ is a integer multiple of $h_\Delta R_\Delta.$
Thus, $\Lambda' = \Lambda$ as soon as $\det (\Lambda') < 2 h^*,$
because $h_\Delta R_\Delta$ is the only integer multiple of itself in
the interval $(h^*,2h^*).$

As described in \cite{JacobsonPhd}, an adaptation of the strategy used
in the self-initialized quadratic sieve (SIQS) factoring algorithm is
used to compute relations.  First, compute the ideal $\mathfrak{a} =
\mathfrak{p}_1^{e_1} \dots \mathfrak{p}_n^{e_n} = (1/d) [ a \Z + (b +
\sqrt{\Delta})/2 \Z]$ with $\mathcal{N}(\mathfrak{a}) = a/d^2.$ Let
$\alpha = (a x + (b + \sqrt{\Delta})/2 y) / d$ with $x,y \in
\mathbb{Z}$ be an arbitrary element in $\mathfrak{a}.$ Then
\begin{equation*}
\mathcal{N}(\alpha) = 
\frac{1}{d^2}
\left( a x + \frac{b + \sqrt{\Delta}}{2} y \right) 
\left( a x + \frac{b - \sqrt{\Delta}}{2} y \right) 
= (a/d^2) (ax^2 + bxy + c y^2)
\end{equation*}
where $c = (b^2 - \Delta)/(4a).$ Because ideal norm is multiplicative,
there exists an ideal $\mathfrak{b}$ with $\mathcal{N}(\mathfrak{b}) =
ax^2 + bxy + cy^2$ such that $(\alpha) = \mathfrak{a} \mathfrak{b}.$
Thus, finding $x$ and $y$ such that $\mathcal{N}(\mathfrak{b})$
factors over the norms of the prime ideals in the factor base yields a
relation.  Such $x$ and $y$ can be found by sieving the polynomial
$\varphi(x,y) = a x^2 + bxy + cy^2,$ and a careful selection of the
ideals $\mathfrak{a}$ yields a generalization of self-initialization,
in which the coefficients of the sieving polynomials and their roots
modulo the prime ideal norms can be computed quickly.  In practice, we
use $\varphi(x,1)$ for sieving, so that the algorithm resembles the
SIQS more closely.  For more details, see \cite{JacobsonPhd} or
\cite{JWPellBook}.

The determinant of the relation lattice $\Lambda'$ is computed in two
stages.  The first step is to compute the determinant of the integer
part of this sublattice by finding a basis in Hermite normal form
(HNF).  Once $\Lambda'$ has full rank, the determinant of this basis
is computed as the product of the diagonal elements in a matrix
representation of the basis vectors.  The group structure is then
computed by finding the Smith normal form of this matrix.  The real
part of $\det (\Lambda'),$ a multiple of the regulator $R_\Delta,$ is
computed by first finding a basis of the kernel of the matrix
consisting of the integer parts of the relations.  Every vector $(k_1,
\dots, k_m) \in \Z^m$ in the kernel corresponds to a multiple of the
regulator computed with $m R_\Delta = k_1 \log |\alpha_1| + \dots +
k_m \log |\alpha_m|.$ The ``real gcd'' of the multiples $m_1 R_\Delta,
\dots, m_n R_\Delta$ computed from each basis vector of the kernel,
defined as $\gcd(m_1,\dots,m_n) R_\Delta,$ is then the real part of
$\det (\Lambda').$ An algorithm of Maurer \cite{maurer} can be used to
compute the real gcd efficiently and with guaranteed numerical
accuracy given explicit representations of the $\alpha_i$ and the
kernel vectors.

As mentioned in the introduction, the correctness of this algorithm
depends on the truth of the Generalized Riemann Hypothesis.  In fact,
the GRH must be invoked in two places.  The first is to compute a
sufficiently accurate approximation $h^*$ of $h_\Delta R_\Delta$ via a
method due to Bach \cite{BachEulerProd}.  Without the GRH, an
exponential number of terms in the Euler product used to approximate
$L(1,\chi_\Delta)$ must be used (see, for example, \cite{Lquadratic}).
The second is to ensure that the factor base generates $Cl_\Delta.$
Without the GRH, an exponential size factor base is required, whereas
by a theorem of Bach \cite{BBounds} the prime ideals of norm less than
$6 \log(\Delta)^2$ suffice.  In practice, an even smaller factor base
is often used, but in that case, the factor base must be verified by
showing that every remaining prime ideal with norm less than Bach's
bound can be factored over the ideals in the factor base.

\section{Practical Improvements}\label{sec:improvements}

In this section, we describe our practical improvements for computing
the class group structure and the regulator of a the real quadratic
field. Some of these improvements, such as the double large prime
variant and structured Gaussian elimination, were used in
\cite{biasse} for the simpler case of imaginary quadratic number
fields. On the other hand, the batch smoothness test and system
solving based methods for computing the regulator had never been
implemented in the context of number fields before.

\subsection{Relation collection}

Improving the relation collection phase allows us to speed up every
other stage of the algorithm. Indeed, the faster the relations are
found, the smaller the factor base can be, thus reducing the
dimensions of the relation matrix and the time taken by the linear
algebra phase. In addition, the verification phase also relies on our
ability to find relations and therefore benefits from improvements to
the relation collection phase. Throughout the rest of the paper, $M$
denotes the relation matrix, the matrix whose rows are the integer
parts of the relations.

\subsubsection{Large prime variants}

The large prime variants were developed in the context of integer
factorization to speed up the relation collection phase in both the
quadratic sieve and the number field sieve. A single large prime
variant was described by Buchmann and D\"ullman \cite{BDDist} for
computing the class group of an imaginary quadratic field, and adapted
to the real case by Jacobson \cite{JacobsonPhd}.  Biasse \cite{biasse} 
described how the double large prime strategy could be
using in the imaginary case, and obtained a significant speed-up.

The idea is to keep relations involving one or two extra primes not in the factor base of
norm less than $B_2\geq B_1$. These relations
thus have the form
%The main idea is the following: We define the ``small primes" to be the prime ideals in the factor base and the small prime bound as the corresponding bound $B_1=B$. Then we define a large prime bound $B_2$. During the relation collection phase we choose not to restrict ourselves to relations only involving primes $\p$ in $\mathcal{B}$ but we also keep relations of the form 
\begin{equation*}
(\alpha)= \p_1^{e_1} \hdots \p_n^{e_n} \p \ \ \text{and}\ \
(\alpha)= \p_1^{e_1} \hdots \p_n^{e_n} \p\p'
\end{equation*} 
for $\p_i$ in $\mathcal{B}$, and for $\p,\p'$ of norm less than
$B_2$. We will refer to these types of partial relations as 1-partial
relations and 2-partial relations, respectively. Keeping partial
relations only involving one large prime is the single large prime
variant, whereas keeping those involving one or two is the double
large prime variant which was first described by Lenstra and Manasse
\cite{Lenstra:2LP}. We do not consider the case of more large primes,
but it is a possibility that has been studied in the context of
factorization \cite{Lenstra:3LP}.

Partial relations may be identified as follows. Let $m$ be the remainder of $\varphi(x,1)$ after the division by all primes $p\leq B_1$, and assume that $B_2 < B_1^2$. If $m=1$ then we have a full relation. If $m\leq B_2$ then we have a 1-partial relation. We can see here that detecting 1-partial relations is almost for free. If we also intend to collect 2-partial relations then we have to consider the following possibilities:
\begin{enumerate}
 \item $m > B_2^2$;
 \item $m$ is prime and $m > B_2$;
 \item $m$ is prime and $m \leq B_2$;
 \item $m$ is composite and $B_1^2 < m \leq B_2^2$.
\end{enumerate}
In Cases 1 and 2 we discard the relation. In Case 3 we have a
1-partial relation, and in Case 4 we have $m=pp'$ where $p =
\mathcal{N}(\p)$ and $p' = \mathcal{N}(\p')$.  Cases~1, 2, and 3 can
be checked very easily, but if none are satisfied we need to factor
$m$ in order to determine whether Case~4 is satisfied.  We used
Milan's implementation of the SQUFOF algorithm \cite{milan} based on
the theoretical work of \cite{squfof} to factor the $m$ values
produced.

Even though we might have to factor the remainder, partial relations are
found much faster than full relations.  However, the dimensions of the
resulting matrix are much larger,
% because the probability that $\mathcal{N}(\mathfrak{b})$ is $B_2$-smooth is much greater than the probability that it is $B_1$-smooth. This improvement in the speed of the relation collection phase comes at a price: The number of columns in the relation matrix is much greater, 
thus preventing us from running the linear algebra phase directly on
the resulting relation matrix. In addition, we have to find many more
relations since we have to produce a full rank matrix. We will see in
\textsection \ref{gauss} how to reduce the dimensions of the relation
matrix using Gaussian elimination techniques.

\subsubsection{Batch smoothness test}

After detecting potential candidates for smooth integers via the SIQS,
one has to certify their smoothness. In \cite{biasse,JacobsonPhd}, this
was done by trial division with the primes in the factor base. The
time taken by trial division can be shortened by using Bernstein's
batch smoothness test \cite{bernstein}, which uses a product tree
structure and modular arithmetic to factor a batch of residues
simultaneously in time $O\left( b(\log b)^2 \log\log b\right)$ where
$b$ is the total number of input bits.

%\begin{algorithm}[H]
%\caption{Batch smoothness test}\label{batch:alg}
%\begin{algorithmic}
%\REQUIRE prime numbers $\left\lbrace p_1,\hdots p_m\right\rbrace $ and positive integers $\left\lbrace  x_1,\hdots x_n\right\rbrace $.
%\ENSURE Array $S$ containing the $\left\lbrace p_1,\hdots p_m\right\rbrace $-smooth part of each $x_i$
%\STATE Compute $z\leftarrow p_1\hdots p_m$ using a product tree.
%\STATE Compute $z \mod x_1,\hdots,z\mod x_n$ using a remainder tree.
%\FOR{$k=1$ to $n$}
%\STATE $e\leftarrow$ smallest integer such that $2^{2^e} \geq x_k$
%\STATE $y_k \leftarrow (z \mod x_k)^{2^e}$
%\STATE $S[k]\leftarrow \text{GCD}(x_k,y_k)$
%\ENDFOR
%\RETURN $S$
%\end{algorithmic}
%\end{algorithm}
%
%\begin{theorem}
%Algorithm \ref{batch:alg} takes time $O\left( b(\log b)^2 \log\log b\right) $ where $b$ is the number of input bits.
%\end{theorem}
%
%\begin{proof}
%See \cite{bernstein}.
%\end{proof}

Instead of testing the smoothness of every potential candidate as soon
as they are discovered, we rather stored them and tested them at the
same time using Bernstein's method as soon their number exceeded a
certain limit. This improvement has an effect that is all the more important
when the time spent in the trial division is long. In our algorithm,
this time mostly depends on the tolerance value $T,$ a parameter used
to control the number of candidates yielded by the sieve for
smoothness testing.
%defined in \textsection \ref{sieving}, as we will see in \textsection
%\ref{numerical}.

\subsection{Structured Gaussian Elimination}\label{gauss}

As mentioned in \textsection \ref{sieving}, in order to determine
whether the computed relations generate the entire relation lattice,
we need to compute the HNF basis of the sublattice they generate.
This can be done by putting the integer components of the relations as
rows in a relation matrix, and computing the HNF.

The first step when using large primes is to compute full relations
from all of the partial relations.  Traditionally, rows were recombined
to give full relations as follows.  In the case of 1-partial
relations, any pair of relations involving the same large prime $\p$
were recombined into a full relation. In the case of 2-partial
relations, Lenstra \cite{Lenstra:2LP} described the construction of a
graph whose vertices were the relations and whose edges linked
vertices having one large prime in common. Finding independent cycles
in this graph allows us to recombine partial relations into full
relations.

In this paper, we instead follow the approach of Cavallar
\cite{cavallar}, developed for the number field sieve, and adapted by
the first author to the computation of ideal class group structures in
imaginary quadratic number fields \cite{biasse}, which uses Gaussian
elimination on columns. The ideas is to eliminate columns using
structured Gaussian strategies until the dimensions of the matrix are
small enough to allow the computation of the HNF with standard
algorithms.

Let us recall a few definitions. First, subtracting two rows is called
\textit{merging}. If two relations corresponding to rows $r_1$ and
$r_2$ share the same prime $\p$ with coefficients $c_1$ and $c_2$
respectively, then multiplying $r_1$ by $c_2$ and $r_2$ by $c_1$ and
merging is called \textit{pivoting}. Finally, finding a sequence of
pivots leading to the elimination of a column of Hamming weight $k$ is
a $k$-way merge.

We aim to reduce the dimensions of the relation matrix by performing
$k$-way merges on the columns of weight $k=1,\hdots,w$ in increasing
order for a certain bound $w$. To limit the growth of the density and
of the size of the coefficients induced by these operations, we used
optimized pivoting strategies. In what follows we describe an
algorithm performing $k$-way merges to minimize the growth of both the
density and the size of the coefficients, thus allowing us to go deeper 
in the elimination process and delay the explosion of the coefficients.

As in \cite{biasse}, we define a cost function $C$ mapping rows onto
the integers. The one used in \cite{biasse} satisfied
\begin{equation}\label{cost:eq}
C(r) = \sum_{1\leq|e_i|\leq Q}1 + c\sum_{|e_j| > Q}1,
\end{equation}
where $c$ and $Q$ are positive numbers, and $r = [e_1,\hdots,e_n]$ is a row
corresponding to $(\alpha)=\prod_i\p_i^{e_i}$. This way, the heaviest
rows are those which have a high density and large coefficients. In our experiments for this
work, we used a different cost function, see \textsection\ref{timings}. Then,
to perform a $k$-way merge on a given column, we construct a complete
graph $\mathcal{G}$ of size $k$ such that
\begin{itemize}
 \item the vertices are the rows $r_i,$ and
 \item every edge linking $r_i$ and $r_j$ has weight $C(r_{ij})$, where $r_{ij}$ is obtained by pivoting $r_i$ and $r_j$.
\end{itemize}
Finding the best sequence of pivots with respect to the chosen cost
function $C$ is equivalent to finding the minimum spanning tree
$\mathcal{T}$ of $\mathcal{G}$, and then recombining every row $r$
with its parent starting with the leaves of $\mathcal{T}$.

Unlike in \cite{biasse}, we need to keep track of the permutations we
apply to the relation matrix, and of the empty columns representing
primes of norm less than $6\log^2\Delta$.  
%Indeed, it will be of use for the kernel computation described in
%\textsection \ref{kernel}.
This will be required for the regulator computation part of the
algorithm described next.

\subsection{Regulator computation}\label{kernel}

As mentioned in \textsection \ref{sieving}, the usual way to compute
the regulator is to find a basis of the kernel of the relation matrix,
compute integer multiples of the regulator from these basis vectors,
and compute their real gcd using Maurer's algorithm \cite{maurer}.
%The usual way of computing the regulator is to find a kernel basis of
%the relation matrix. Then, any power-product of the $(\alpha_i)$
%represented by the rows with exponents corresponding to entries of a
%kernel vector is a unit. A multiple of the regulator can be computed
%with the function \texttt{Regulator multiple} which description can be
%found in \cite{maurer}. It takes the generetors of the relations and a
%set of kernel vectors in input and outputs a multiple of $R$. We
%obtain the regulator this way provided the elements
%$$v_j := u^j_1\log|\alpha_1| + \hdots + u^j_n\log|\alpha_n|$$ can be
%recombined to form $R$, where the vectors $(u^j_i)_{i\leq n}$ for
%$j=0,\hdots \dim(\ker M)$ are a kernel basis of the relation
%matrix. 
If $\det \Lambda' > 2h^*,$ then either the class number or regulator
computed is too large, and we need to find extra relations
corresponding to new generators, and new kernel vectors involving them.

In this section, we describe a way of taking advantage of the large
number of generators involved in the different partial
relations. Indeed, the dimensions of the relation matrix before the
Gaussian elimination stage is much larger than in the base scenario
and thus involves more generators. Consequently, given a set of $k\leq
\dim (\ker M)$ kernel vectors $(u^j_1,\hdots,u^j_n)_{j\leq k}$, the
probability that the corresponding elements
\begin{equation*}
v_j := u^j_1\log|\alpha_1| + \hdots + u^j_n\log|\alpha_n| \enspace,
\end{equation*}
where $\alpha_i$ is the generator of the $i$-th relation, can be
recombined into $R$ is much larger. On the other hand, the dimensions
of the matrix prevents us from running a kernel computation directly
after the relation collection phase. 
%We could be tempted to do it on the small matrix whose dimension is
%lower and then use the permutation to turn them into kernel vectors of
%the original matrix, bu t we would systematically obtain multiples of
%the regulator since the non-zero entries of such vectors correspond to
%a small subset of the generators.
Thus, rather than attempting to compute the kernel, we use a method
similar to that of Vollmer \cite{Vollmer-regulator} based on solving
linear systems.

The first step of our algorithm consists of putting the matrix in a
pseudo-lower triangular form using a permutation obtained during the
Gaussian elimination phase. Indeed, as part of this computation we
obtain a unimodular matrix $U\in\Z^{n\times n}$ such that
\[ UM = \left( 
   \begin{BMAT}(@)[2pt,3cm,3.5cm]{c.c}{c.c}
   \begin{BMAT}[2pt,1cm,1.5cm]{c}{c} 
	A
	\end{BMAT} & 
	\begin{BMAT}[2pt,2cm,1.5cm]{c}{c} 
	(0)
	\end{BMAT} \\
   \begin{BMAT}[2pt,1cm,2.5cm]{c}{c} 
	(*)
   \end{BMAT} & 
   \begin{BMAT}[2pt,2cm,2.5cm]{ccc}{ccc}
1&       & (0)        \\
 &   \ddots &       \\
(*)       &  & 1 
  \end{BMAT} 
\end{BMAT}
   \right).  \]
Thus, solving a linear system of the form $xM = b$ for a vector
$b\in\Z^m$ boils down to solving a system of the form $x'A = b'$, then
doing a trivial descent through the diagonal entries which equal 1
and finally permuting back the coefficients using~$U$. To solve the
small linear systems, we used the algorithm \texttt{certSolveRedLong}
from the IML library \cite{iml}. It takes a single precision dense
representation of $A$ and returns an LLL-reduced solution.

Once $M$ is in pseudo-lower triangular form, we draw a set of
relations $r_1,\hdots r_d$ which are not already rows of $M$, and
for each $r_i$, $i\leq d$, we solve the system $x_i A = r_i.$ We then
augment $M$ with the rows $r_i$ for $i\leq d$ and the vectors $x_i$ with~$d$ 
extra coordinates, which are all set to zero except for the $i$-th
which is set to~$-1$.
\[ M' := \left( 
   \begin{BMAT}(@)[2pt,1cm,2cm]{c}{c.c}
   \begin{BMAT}(e){c}{c}
M
  \end{BMAT} \\
\begin{BMAT}[2pt,1cm,0.5cm]{c}{c} 
	r_i
\end{BMAT}
\end{BMAT}
   \right)  \ \ 
x_i' := \left(
\begin{BMAT}(@)[2pt,3cm,0.5cm]{c.c}{c}
\begin{BMAT}(e){c}{c}
x_i
  \end{BMAT} & 
\begin{BMAT}[2pt,1cm,0.5cm]{ccc}{c}
0\hdots 0 & -1 & 0\hdots 0
\end{BMAT} 
\end{BMAT}\right). 
\]
We clearly have $x_i' M' = 0$ for $i\leq d,$ and the $x_i'$ can be
used to find a multiple of $R_\Delta$ as described in \textsection
\ref{sieving}.
%In addition, the generators involved in those kernel vectors are more
%randomized, thus allowing us to use fewer vectors. To conclude, let us
%recall the main steps of this algorithm:
%
%\begin{algorithm}[H]
%\caption{Regulator computation}\label{reg:alg}
%\begin{algorithmic}
%\REQUIRE $M$, $U$, $d$.
%\ENSURE $R$.
%\STATE Find $d$ extra relations.$(r_i)_{i\leq d}$.
%\STATE Solve the $d$ systems $X_iM = r_i$.
%\STATE Construct $(X'_i)_{i\leq d}$.
%\STATE Give the $(X'_i)_{i\leq d}$ to \texttt{Regulator multiple} along with the corresponding generators,
%\RETURN $R$
%\end{algorithmic}
%\end{algorithm}

\section{Numerical results}\label{numerical}

In this section, we give numerical results showing the impact of our
improvements. For each timing, we specify the
architecture used. All the timings were obtained with our code in C++
based on the libraries GMP \cite{gmp}, NTL \cite{NTL}, IML \cite{iml}
and Linbox \cite{linbox}. All timings are in
CPU seconds.

\subsection{Comparative timings}\label{timings}

The state of the art concerning class group and regulator computation
was established in \cite{JacobsonPhd}, where all the timings were
obtained with the SPARCStation II architecture. In addition, most of the
code used at the time is unavailable now, including the HNF
computation algorithm. Thus, providing a meaningful comparison between our
methods and those of \cite{JacobsonPhd} is difficult. We chose to
implement the HNF computation algorithm in a way that resembles the one
of \cite{JacobsonPhd}, but takes advantage of the libraries
available today for computing the determinant and the modular HNF. We used 
this implementation in each different scenario. The
relation collection phase is easier to compare, since our method relies on SIQS.

In the following, we will refer to the base case as the strategy consisting 
of finding the relation matrix without using the large prime variants or the 
smoothness batch test, and calculating the regulator by computing its kernel 
with the algorithm \texttt{nullspaceLong} from IML library. It differs from 
the 0 large prime case (0LP) where we use the algorithm described in \textsection
\ref{kernel} for computing the regulator, along with a relation 
collection phase that does not use large primes. We also denote the 1 large prime scenario 
by 1LP, the 2 large primes by 2LP and 2LP Batch when using batch smoothness test.

\subsubsection{Relation collection phase}

In Table \ref{tab:sieve}, we give the time taken to collect all
necessary relations. Without large primes, we collected
$|\mathcal{B}|+100$ relations, whereas when we allow large primes we
need to collect enough relations to ensure that the number of rows is
larger than the number of non-empty columns.  We used a 2.4 GHz
Opteron with 16GB of memory and took $\Delta = 4(10^n+3)$ with~$40\leq
n\leq 70$. For each discriminant, we used the optimal parameters given
in \cite{JacobsonPhd}, including the size of the factor base, even if
we tend to reduce this parameter when optimizing the overall time.
The only parameter we modified is the tolerance value for the SIQS, as
a higher tolerance value is required for the large prime
variations. In each case we took $B_2 = 12 B_1$.  It is shown in
\cite{biasse} that the ratio $B_2/B_1$ does not have an important
impact on the sieving time.

\begin{table}[!ht]
\caption{Comparative table of the relation collection time}
\label{tab:sieve}
%\small 
\begin{center}
 \begin{tabular}{|r|r|r|r|r|}
\hline
\multicolumn{1}{|c|}{$n$} & 
\multicolumn{1}{c|}{0LP} & 
\multicolumn{1}{c|}{1LP} & 
\multicolumn{1}{c|}{2LP} &
\multicolumn{1}{c|}{2LP Batch} \\
\hline
40 &     0.83 &    0.48 &    0.63 &   0.90 \\
45 &     6.70 &    3.10 &    2.70 &   2.20 \\
50 &    23.00 &    9.50 &    9.20 &   6.10 \\
55 &    56.00 &   26.00 &   23.00 &  15.00 \\
60 &   202.00 &   86.00 &   69.00 &  41.00 \\
65 &  1195.00 &  513.00 &  354.00 & 227.00 \\
70 &  4653.00 & 1906.00 & 1049.00 & 834.00 \\
\hline
\end{tabular}
\end{center}
\end{table}

The timings in Table~\ref{tab:sieve} correspond to the optimal value
of the tolerance value in each case, found by trying values between 1.7 and 4, and keeping the optimum for each scenario. 
For 0LP, the optimal value is between 1.7 and 2.3
whereas it is around 2.3 for 1LP, 2.8 for 2LP
and 3.0 for 2LP Batch. The latter
case has a higher optimal tolerance value because using the batch
smoothness test allows one to spend more time factoring the
residues. When using Bernstein's smoothness test, we took batches of
100 residues. In our experiments, this value did not seem to have an important
effect on the relation collection time. We observe in Table \ref{tab:sieve} that the use of
the large prime variants has a strong impact on the relation
collection phase, and that using the smoothness batch test strategy
yields an additional speed-up of approximately 20\% over the double
large prime strategy.

\subsubsection{Structured Gaussian elimination}

Structured Gaussian elimination allows us to reduce the time
taken by the linear algebra phase by reducing the dimensions of the
relation matrix. Our method minimizes the growth of the density and of
the size of the coefficients. To illustrate the impact of the
algorithm described in \textsection \ref{gauss}, we monitor in Table
\ref{crunch:tab} the evolution of the dimensions of the matrix, the
average Hamming weight of its rows, the extremal values of its
coefficients and the time taken for computing its HNF in the case of a
relation matrix corresponding to $\Delta = 4(10^{60} +3)$. We keep
track of these values after all $i$-way merges for some values of $i$
between 5 and 170. The original dimensions of the matrix are $2000\times
1700$, and the timings are obtained on a 2.4 Ghz Opteron with 32GB of
memory.

In \cite{biasse}, the first author regularly deleted the rows having
the largest coefficients. To do this, we need to create more rows than
in the base case. To provide a fair comparison between the two
strategies, we used the same relation matrix resulting from a relation
collection phase without large primes, and with as few
rows as was required to use the same algorithm as in \cite{JacobsonPhd}. We
therefore had to drop the regular row deletion. We also tuned
the cost function to compensate for the resulting growth of the
coefficients, using
\begin{equation*}
C(r) = \sum_{1\leq|e_i|\leq 8}1 + 100\sum_{|e_j| > 8}|e_j| \enspace,
\end{equation*}
instead of \eqref{cost:eq}. 

The HNF computation consists of taking the GCD of the determinants of
two different submatrices of the matrix after elimination using
Linbox, and using the modular HNF of NTL with this value. Indeed, this
GCD (which is likely to be relatively small) is a multiple of
$h_\Delta$. This method, combined with an elimination strategy due to
Havas \cite{HMdiag}, was used in \cite{JacobsonPhd} and implemented in
LiDIA \cite{LiDIA97}.  As this implementation is no longer available,
we instead refer to the timings of our code, which has the advantage
of using the best linear algebra libraries available today.

\begin{table}[!ht]
\caption{Comparative table of elimination strategies}
\label{crunch:tab}
\small 
\begin{center}
 \begin{tabular}{|r|r|r|r|r|r|r|}
\hline
  \multicolumn{7}{|c|}{\textbf{Naive Gauss}} \\
\hline
\multicolumn{1}{|c|}{$i$} & 
\multicolumn{1}{|c|}{Row Nb} & 
\multicolumn{1}{|c|}{Col Nb} & 
\multicolumn{1}{|c|}{Average weight} & 
\multicolumn{1}{|c|}{max coeff} & 
\multicolumn{1}{|c|}{min coeff} & 
\multicolumn{1}{|c|}{HNF time} \\
\hline
  5 & 1189 & 1067 &  27.9 &  14 &  -17 & 357.9 \\
 10 &  921 &  799 &  49.3 &  22 &  -19 & 184.8 \\
 30 &  757 &  635 & 112.7 &  51 &  -50 & 106.6 \\
 50 &  718 &  596 & 160.1 &  81 &  -91 &  93.7 \\
 70 &  699 &  577 & 186.3 & 116 & -104 &  85.6 \\
 90 &  684 &  562 & 205.5 & 137 &  -90 &  79.0 \\
125 &  664 &  542 & 249.0 & 140 & -146 &  73.8 \\
160 &  655 &  533 & 282.4 & 167 & -155 &  72.0 \\
170 &  654 &  532 & 286.4 & 167 & -155 & 222.4 \\
\hline
\multicolumn{7}{|c|}{\textbf{With dedicated elimination strategy}} \\
\hline
\multicolumn{1}{|c|}{$i$} & 
\multicolumn{1}{|c|}{Row Nb} & 
\multicolumn{1}{|c|}{Col Nb} & 
\multicolumn{1}{|c|}{Average weight} & 
\multicolumn{1}{|c|}{max coeff} & 
\multicolumn{1}{|c|}{min coeff} & 
\multicolumn{1}{|c|}{HNF time} \\
\hline
  5 & 1200 & 1078 &  26.8 &       13 &       -12 & 368.0 \\
 10 &  928 &  806 &  42.6 &       20 &       -15 & 187.2 \\
 30 &  746 &  624 &  82.5 &       33 &       -27 & 100.8 \\
 50 &  702 &  580 & 107.6 &       64 &       -37 &  84.3 \\
 70 &  672 &  550 & 136.6 &      304 &      -676 &  73.4 \\
 90 &  656 &  534 & 157.6 &     1278 &     -1088 &  67.5 \\
125 &  637 &  515 & 187.1 &     3360 &     -2942 &  63.4 \\
160 &  619 &  497 & 214.6 &     5324 &     -3560 &  56.9 \\
170 &  615 &  493 & 247.1 & 36761280 & -22009088 & 192.6 \\
\hline
 \end{tabular}
\end{center}
\end{table}
\normalsize

Table \ref{crunch:tab} shows that the use of our elimination strategy
leads to a matrix with smaller dimensions (493 rows with our method,
533 with the naive elimination) and lower density (the average weight
of its rows is of 214 with our method and 282 with the naive
elimination). These differences result in an improvement of the time
taken by the HNF computation: 56.9 seconds with our method against
72.0 seconds with the naive Gaussian elimination. The regular
cancellation of the rows having the largest coefficients over the
course of the algorithm would delay the explosion of the coefficient
size, but require more rows for the original matrix. This brutal
increase in the size of the extremal values of the matrix can be seen
in Table \ref{crunch:tab}.  At this point these higher values
propagate during pivoting operations, and any further column
elimination becomes counter-productive.

\subsubsection{Factor base verification}

The improvements in the relation collection phase have an impact on
the factor base verification. The impact of the smoothness batch test is straightforward,
whereas the large prime variants act in a more subtle way. Indeed,
we create many more relations when using the large prime variants, and
the relations created involve primes of larger norm. Therefore, a
given prime not in $\mathcal{B}$ of norm less than $6\log^2\Delta$ is more likely to appear in a relation, and thus not to
need to be verified. Table \ref{tab:verif} shows the impact of 
the large prime variants and of the batch smoothness test on the
verification time. We used a 2.4 GHz Opteron with 16GB of memory. We
considered discriminants of the form $\Delta = 4(10^n+3)$ for $n$
between 40 and 70, and we chose in every case the factor base giving
the best results for the base scenario.

\begin{table}[!ht]
\caption{Comparative table of the factor base verification time}
\label{tab:verif}
%\small 
\begin{center}
 \begin{tabular}{|r|r|r|r|r|}
\hline
\multicolumn{1}{|c|}{$n$} & 
\multicolumn{1}{c|}{0LP} & 
\multicolumn{1}{c|}{1LP} & 
\multicolumn{1}{c|}{2LP} &
\multicolumn{1}{c|}{2LP Batch} \\
\hline
40 &   17.0 &   11.0 &   11.0 &    6.2 \\
45 &   77.0 &   44.0 &   30.0 &   18.0 \\
50 &  147.0 &   85.0 &   52.0 &   43.0 \\
55 &  308.0 &  167.0 &  134.0 &  110.0 \\
60 &  826.0 &  225.0 &  282.0 &  274.0 \\
65 & 8176.0 & 1606.0 & 1760.0 & 1689.0 \\
70 & 9639.0 & 4133.0 & 5777.0 & 2706.0 \\
\hline
\end{tabular}
\end{center}
\end{table}

\subsubsection{Regulator computation}

Our method for computing the regulator avoids computing the relation
matrix kernel. Instead, we need to solve a few linear systems
involving the matrix resulting from the Gaussian elimination. To
illustrate the impact of this algorithm, we used the relation matrix
obtained in the base case for discriminants of the form $4(10^n +3 )$
for $n$ between 40 and 70. The timings are obtained on a 2.4GHz
Opteron with 16GB of memory.

\begin{table}[!ht]
\caption{Comparative table of regulator computation time}
\label{tab:reg}
%\small 
\begin{center}
 \begin{tabular}{|r|r|r|}
\hline
\multicolumn{1}{|c|}{$n$} & 
\multicolumn{1}{c|}{Kernel Computation} & 
\multicolumn{1}{c|}{System Solving} \\
\hline
40 &    15.0 &   6.2 \\
45 &    18.0 &   8.3 \\
50 &    38.0 &  20.0 \\
55 &   257.0 &  49.0 \\
60 &   286.0 & 103.0 \\
65 &  5009.0 & 336.0 \\
70 & 10030.0 & 643.0 \\
\hline
\end{tabular}
\end{center}
\end{table}

In Table \ref{tab:reg}, the timings corresponding to our system
solving approach are taken with seven kernel vectors.
However, in most cases only two or three vectors are required to compute
the regulator. As most of the time taken by our approach is spent on
system solving, we see that computing fewer kernel vectors would
result in an improvement of the timings, at the risk of obtaining a
multiple of the regulator.

\subsubsection{Overall time}

We have studied the individual impact of our improvements on each
stage of the algorithm. We now present their effect on the overall
time taken by the algorithm, including the factor base verification
time, for discriminants of the form $\Delta = 4(10^n+3)$ with $40\leq
n\leq 70$ on a 2.4 GHz Opteron with 16GB of memory. We used the same
parameters as in \cite{JacobsonPhd}, except for the tolerance and the
size of the factor base. We notice in Table \ref{overall:tab} that the
optimal size of the factor base is smaller when we use improvements
for the sieving phase. For example the optimal size for the double
large prime variant is half the one of the base case scenario. This
results in an improvement in the HNF and regulator computation whereas
the relation collection time can remain unchanged, or even
increase. The tolerance value we chose varies only with the strategy,
but not with the size of the discriminant. We chose 2.0 for the base
case and 0LP whereas we set it to 2.3 for 1LP, 2.8 for 2LP and 3.0 for
2LP Batch. We eliminated columns of weight up to $w = 150$ since Table~\ref{crunch:tab}
indicates that further elimination is counter-productive.

\begin{table}[!ht]
\caption{Effect on the overall time}
\label{overall:tab}
\small 
\begin{center}
 \begin{tabular}{|r|r|r|r|r|r|r|r|r|}
\hline
\multicolumn{1}{|c|}{$n$} & 
\multicolumn{1}{c|}{strategy} & 
\multicolumn{1}{c|}{$|\mathcal{B}|$} &
\multicolumn{1}{c|}{relations} & 
\multicolumn{1}{c|}{elimination} & 
\multicolumn{1}{c|}{HNF} & 
\multicolumn{1}{c|}{regulator} &
\multicolumn{1}{c|}{verification} & 
\multicolumn{1}{c|}{total} \\
\hline
 \multirow{4}{*}{40} &  base & 400& 0.8 & 0.1 & 3.2 & 14.6 & 16.8 & 35.6 \\
  & 0LP & 400& 0.7 & 0.1 & 2.2 & 6.0 & 16.6  & 25.7 \\
  & 1LP & 300& 0.8 & 0.2 & 2.5 & 6.4 & 13.1 & 23.1 \\
  & 2LP & 250& 1.7 & 0.3 & 4.8 & 8.7 & 18.0  & 33.3 \\
  & 2LP Batch & 250&  0.5 & 0.2 & 3.6 & 6.7 & 4.4 & 15.5 \\
\hline
\multirow{4}{*}{45} &  base & 500&  6.7 & 0.1 & 5.1 & 18.0 & 77.0 & 107.0 \\
  & 0LP & 500&  5.9 & 0.2 & 4.9 & 10.0 & 85.0 &  106.0 \\
  & 1LP & 400&  4.0 & 0.4 & 6.0 & 11.0 & 50.0 & 71.0 \\
  & 2LP & 350&  3.8 & 0.5 & 12.0 & 17.0 & 36.0 & 69.0 \\
  & 2LP Batch  & 350&  2.6 & 1.1 & 9.0 & 14.0 & 30.0  & 57.0\\
\hline
\multirow{4}{*}{50} &  base & 750&  23.0 &  0.3 & 16.0 & 38.0 & 147.0 & 224.0  \\
  & 0LP & 700& 21.0 & 0.4 & 15.0 & 20.0 & 147.0 & 203.0 \\
  & 1LP & 450&  20.0 & 0.4 & 10.0 & 17.0 & 108.0 & 155.0 \\
  & 2LP & 400&  14.0 & 0.8 & 22.0 & 23.0 & 74.0 & 133.0 \\
  & 2LP Batch & 400& 10.0 & 0.6 & 21.0 & 25.0 & 62.0 & 119.0 \\

\hline
\multirow{4}{*}{55} &  base & 1200&  129.0 & 1.9 & 60.0 & 257.0 & 308.0 & 756.0 \\
  & 0LP & 1300& 47.0 & 0.7 & 52.0 & 49.0 & 265.0 & 414.0 \\
  & 1LP & 650& 61.0 & 0.7 & 28.0 & 33.0 & 255.0 & 378.0 \\
  & 2LP & 550& 40.0 & 1.1 & 48.0 & 48.0 & 177.0 & 313.0  \\
  & 2LP Batch & 550& 34.0 & 1.0 & 47.0 & 48.0 & 141.0 & 271.0 \\

\hline
\multirow{4}{*}{60} &  base & 1700&  322.0 & 2.9 & 95.0 & 286.0 & 830.0 & 1535.0  \\
  & 0LP & 1700&  187.0 & 1.3 & 106.0 & 103.0 & 846.0 & 1244.0 \\
  & 1LP & 750&  309.0 & 1.0 & 45.0 & 64.0 & 865.0 &  1284.0 \\
  & 2LP& 700 &  143.0 & 2.1 & 152.0 & 137.0 & 365.0 & 799.0 \\
  & 2LP Batch & 700& 142.0 & 1.8 & 103.0 & 100.0 & 309.0 & 655.0 \\
\hline
\multirow{4}{*}{65} &  base & 2700& 10757.0 & 12.0 & 652.0 & 5009.0 & 8176.0& 24607.0  \\
  & 0LP & 2700& 1225.0 & 2.8 & 489.0 & 336.0 & 3676.0 & 5730.0 \\
  & 1LP & 1900 & 1003.0 & 15.0 & 318.0 & 262.0 & 2984.0 & 4583.0 \\
  & 2LP & 1200& 753.0 & 4.7 & 525.0 & 398.0 & 1943.0 & 3624.0 \\
  & 2LP Batch & 1000& 1030.0 & 35.0 & 199.0 & 219.0 & 1642.0 & 3125.0 \\
\hline
\multirow{4}{*}{70} &  base& 3700 & 17255.0 & 24.0 & 1869.0 & 10031.0 & 9639.0 & 38818.0 \\
  & 0LP & 3600&  4934.0 & 19.0 & 1028.0 & 644.0 & 9967.0 &  16591.0 \\
  & 1LP  & 2500&  3066.0 & 17.0 & 845.0 & 646.0 & 9005.0 & 13579.0\\
  & 2LP & 1700&  2414.0 & 27.0 & 2054.0 & 1295.0 & 4590.0 & 10379.0 \\
  & 2LP Batch & 1700& 2588.0 & 20.0 & 1372.0 & 934.0 & 5078.0 & 9991.0 \\
\hline
 \end{tabular}
\end{center}
\end{table}
\normalsize

Table \ref{overall:tab} shows that there is an overall speed-up of of
a factor of 2 for the smallest discriminants and 4 for the
largest. The base case with the largest discriminants suffers from the
necessity of finding some relations in a more randomized way. This
ensures that we can get full rank submatrices of the relation matrix
after the Gaussian elimination to compute a small multiple of
$h_{\Delta}$. Matrices produced using the large prime variants do not
need this extra step, even with the largest discriminants. This
naturally affects the sieving time, since we cannot use SIQS for that
purpose, but also affects phases relying on linear algebra. Indeed,
elimination produces a matrix with larger entries and dimensions.

\subsection{Large example}

The improvements we described allow us to compute class groups and
regulators of real number fields with larger discriminants than was
previously possible. The key is to parallelize the relation collection
and verification phase, while the linear algebra has to be performed
the usual way. These methods were successfully used in \cite{biasse}
to compute the class group structure of an imaginary quadratic field
with a 110-digit discriminant.  We used a cluster with 260 2.4GHz Xeon
cores to compute a relation matrix corresponding to the discriminant
$\Delta_{110} := 4(10^{110}+3)$ in 4 days. We allowed two large primes, used a tolerance value of 3.0, 
tested batches of 100 residues, took $w = 250$ and set
$|\mathcal{B}| = 13000$ . Then, we used three 2.4 GHz
Opterons with 32GB of memory each to compute determinants of full-rank
submatrices of the relation matrix after the Gaussian elimination in 1
day, and one 2.4GHz Opteron to compute the HNF modulo the GCD of these
determinants in 3 days. We had to find 4018 extra relations during the verification phase that
took 4 days on 96 2.4GHz Xeon cores. We thus obtained that
\begin{equation}\label{d110}
Cl_{\Delta_{110}} \cong \Z/12\Z \times \Z/2\Z \enspace,
\end{equation}
and the corresponding regulator is 
\small
\begin{equation*}
R_{\Delta_{110}} \approx
70795074091059722608293227655184666748799878533480399.6730200233
\enspace.
\end{equation*}
\normalsize We estimate that it would take two weeks (4000 relations
per day) to complete the relation collection for $\Delta_{120}$ with
the same factor base as $\Delta_{110}$, thus requiring a similar time
for the linear algebra.

\section{Conclusions}

Recently, our work has been extended to the problems of principal
ideal testing and solving the discrete logarithm problem in the ideal
class group \cite{BJSestimates}.  The double large prime variant and
improvements to relation generation translated directly to
improvements in this context.  However, HNF computations are not
required for this problem, and linear system solving over $\mathbb{Z}$
can be used instead.  The numerical results were used to give
estimates for discriminant sizes that offer equivalent security to
recommended sizes of RSA moduli.

Some possibilities for further improvements remain to be investigated.
For example, a lattice sieving strategy could be used to sieve
$\varphi(x,y)$ instead of $\varphi(x,1).$ Factor refinement and
coprime factorization techniques may be a useful alternative to
Bernstein's batch smoothness test.  Multiple large primes have been
successfully used for integer factorization and could also be tried in
our context.

There is also still room for improvement to the linear algebra
components.  For example, a HNF algorithm that exploits the natural
sparseness of the relation matrix, perhaps as a black-box algorithm,
would be useful. If such an algorithm were available, we could
reconsider using Gaussian elimination techniques since they induce a
densification of the matrix. We could also study the effect of other
dense HNF algorithms in existing linear algebra packages such as KASH,
Pari, Sage and especially MAGMA which seems to have the most efficient
HNF algorithm for our types of matrices. In that case, we would need
the elimination phase regardless of how these algorithms are affected
by the density and the size of the coefficients of the matrix. Indeed,
we cannot afford manipulating a dense representation of the matrix
before the Gaussian elimination phase.

\providecommand{\bysame}{\leavevmode\hbox to3em{\hrulefill}\thinspace}
\providecommand{\MR}{\relax\ifhmode\unskip\space\fi MR }
\providecommand{\MRhref}[2]{%
  \href{http://www.ams.org/mathscinet-getitem?mr=#1}{#2}
}
\providecommand{\href}[2]{#2}


\begin{thebibliography}{10}

\bibitem{BBounds}
E.~Bach, \emph{Explicit bounds for primality testing and related problems},
  Math. Comp. \textbf{55} (1990), no.~191, 355--380.

\bibitem{BachEulerProd}
\bysame, \emph{Improved approximations for {E}uler products}, Number Theory:
  CMS Proc., vol.~15, Amer. Math. Soc., Providence, RI, 1995, pp.~13--28.

\bibitem{bernstein}
D.~Bernstein, \emph{How to find smooth parts of integers}, submited to
  \textit{Mathematics of Computation}.

\bibitem{biasse}
J-F. Biasse, \emph{Improvements in the computation of ideal class groups of
  imaginary quadratic number fields}, to appear in \textit{Advances in
  Mathematics of Communications}, 2010.

\bibitem{BJSestimates}
J-F. Biasse, M. J. Jacobson, Jr., A. K. Silvester, \emph{Security
  estimates for quadratic field based cryptosystems}, to appear in ACISP 2010.

\bibitem{BSub}
J.~Buchmann, \emph{A subexponential algorithm for the determination of class
  groups and regulators of algebraic number fields}, S\'{e}minaire de
  Th\'{e}orie des Nombres (Paris), 1988-89, pp.~27--41.

\bibitem{BDDist}
J.~Buchmann and S.~D{\"u}llmann, \emph{Distributed class group computation},
  Festschrift aus {A}nla{\ss} des sechzigsten {G}eburtstages von {H}errn
  {P}rof. {D}r. {G}. {H}otz, Universit{\"a}t des Saarlandes, 1991, and Teubner,
  Stuttgart, 1992, pp.~69--79.

\bibitem{cavallar}
S.~Cavallar, \emph{Strategies in filtering in the number field sieve}, ANTS-IV:
  Proceedings of the 4th International Symposium on Algorithmic Number Theory,
  Lecture Note in Computer Science, vol. 1838, Springer-Verlag, 2000,
  pp.~209--232.

\bibitem{iml}
Z.~Chen, A.~Storjohann, and C.~Fletcher, \emph{{IML: Integer Matrix Library}},
  \mbox{Software}, 2010, see \url{http://www.cs.uwaterloo.ca/~astorjoh/iml.html}.

\bibitem{Lenstra:3LP}
B.~Dodson, P.~C. Leyland, A.~K. Lenstra, A.~Muffett, and S.~Wagstaff,
  \emph{M{P}{Q}{S} with three large primes}, ANTS-V: Proceedings of the 5th
  International Symposium on Algorithmic Number Theory, Lecture Note in
  Computer Science, vol. 2369, Springer-Verlag, 2002, pp.~446--460.

\bibitem{gmp}
GMP, \emph{The {GNU} multiple precision bignum library}, \mbox{Software}, 2010,
  see \url{http://gmp\-lib.org/}.

\bibitem{squfof}
J.~E. Gower and S.~Wagstaff, \emph{Square form factorization}, Mathematics of
  Computation \textbf{77} (2008), 551--588.

\bibitem{HMdiag}
G.~Havas and B.S. Majewski, \emph{Integer matrix diagonalization}, Journal of
  Symbolic Computing \textbf{24} (1997), 399--408.

\bibitem{JacobsonPhd}
M.~J. Jacobson, Jr., \emph{Subexponential class group computation in quadratic
  orders}, Ph.D. thesis, Technische Universität Darmstadt, Darmstadt, Germany,
  1999.

\bibitem{JSWkeyex}
M.~J. Jacobson, Jr., R.~Scheidler, and H.~C. Williams, \emph{The efficiency and
  security of a real quadratic field based key exchange protocol}, Public-Key
  Cryptography and Computational Number Theory (Warsaw, Poland), de Gruyter,
  2001, pp.~89--112.

\bibitem{JWPellBook}
M.~J. Jacobson, Jr. and H.~C. Williams, \emph{Solving the {P}ell equation}, CMS
  Books in Mathematics, Springer-Verlag, 2009, ISBN 978-0-387-84922-5.

\bibitem{Lenstra:2LP}
A.~K. Lenstra and M.~S. Manasse, \emph{Factoring with two large primes
  (extended abstract)}, Advances in Cryptology - EUROCRYPT '90, Lecture Note in
  Computer Science, vol. 473, Springer-Verlag, 1991, pp.~72--82.

\bibitem{LiDIA97}
LiDIA Group, \emph{{LiDIA}: a {\tt c++} library for computational number
  theory}, Software, Technische Universit{\"a}t Darmstadt, Germany,
  1997, see \url{http://www.informatik.tu-darm\-stadt.de/TI/LiDIA}.

\bibitem{linbox}
LinBox, \emph{Project {LinBox}: {E}xact computational linear algebra},
  \mbox{Software}, 2010, see \url{http://www.linalg.org/}.

\bibitem{Lquadratic}
S.~Louboutin, \emph{Computation of class numbers of quadratic number fields},
  Math. Comp. \textbf{71} (2002), no.~240, 1735--1743.

\bibitem{maurer}
M.~Maurer, \emph{Regulator approximation and fundamental unit computation for
  real quadratic orders}, Ph.D. thesis, Technische Universität Darmstadt,
  Darmstadt, Germany, 1999.

\bibitem{milan}
J.~Milan, \emph{Tifa}, \mbox{Software}, 2010,
  http://www.lix.polytechnique.fr/Labo/Jerome.\-Milan/tifa/tifa.xhtml.

\bibitem{NTL}
V.~Shoup, \emph{{NTL: A Library for doing Number Theory}}, \mbox{Software},
  2010, \url{http://www.\-shoup.net/ntl}.

\bibitem{Vollmer-regulator}
U.~Vollmer, \emph{An accelerated {B}uchmann algorithm for regulator computation
  in real quadratic fields}, Algorithmic Number Theory --- ANTS-V, Lecture
  Notes in Computer Science, vol. 2369, 2002, pp.~148--162.

\end{thebibliography}
\end{document}